\documentclass[12pt]{amsart}
\usepackage{amsmath,amssymb,latexsym}
\usepackage[francais]{babel}
\usepackage[latin1]{inputenc}

\def\p{\varphi}

\def\acl{\mathrm{acl}}

\def\QQ{\mathbb Q}

\def\In{\mathrm{In}}
\def\Pin{\Psi}
\def\aIn{\widetilde{\mathrm{In}}}
\def\aPin{\widetilde\Psi}
\def\F{\mathcal F}

\def\bdd{\mathrm{bdd}}

\def\acl{\mathrm{acl}}
\def\tp{\mathrm{tp}}

\def\stp{\mathrm{stp}}

\def\Ind#1#2{#1\setbox0=\hbox{$#1x$}\kern\wd0\hbox to 0pt{\hss$#1\mid$\hss}
\lower.9\ht0\hbox to 0pt{\hss$#1\smile$\hss}\kern\wd0}
\def\ind{\mathop{\mathpalette\Ind{}}}
\def\Notind#1#2{#1\setbox0=\hbox{$#1x$}\kern\wd0\hbox to 0pt{\mathchardef
\nn="3236\hss$#1\nn$\kern1.4\wd0\hss}\hbox to 0pt{\hss$#1\mid$\hss}\lower.9\ht0
\hbox to 0pt{\hss$#1\smile$\hss}\kern\wd0}

\theoremstyle{plain}
\newtheorem{theorem}{Théorème}[section]
\newtheorem{proposition}[theorem]{Proposition}
\newtheorem{fact}[theorem]{Fait}
\newtheorem{lemma}[theorem]{Lemme}
\newtheorem{corollary}[theorem]{Corollaire}

\theoremstyle{definition}
\newtheorem{definition}[theorem]{Définition}
\newtheorem{remark}[theorem]{Remarque}
\newtheorem{vermutung}[theorem]{Conjecture}
\newtheorem{beispiel}[theorem]{Exemple}
\def\defn{\begin{definition}}
\def\edefn{\end{definition}}
\def\satz{\begin{theorem}}
\def\esatz{\end{theorem}}
\def\tats{\begin{fact}}
\def\etats{\end{fact}}
\def\kor{\begin{corollary}}
\def\ekor{\end{corollary}}
\def\lmm{\begin{lemma}}
\def\elmm{\end{lemma}}
\def\bem{\begin{remark}}
\def\ebem{\end{remark}}
\def\bew{\par\noindent{\em Démonstration : }}
\def\satzli{\begin{proposition}}
\def\esatzli{\end{proposition}}
\def\bsp{\begin{beispiel}}
\def\ebsp{\end{beispiel}}
\def\verm{\begin{vermutung}}
\def\everm{\end{vermutung}}

\begin{document}
\title{Quelques réflexions inévitables}
\author{Frank Wagner}
\address{Universit\'e de Lyon; Universit\'e Lyon 1; CNRS; Institut Camille Jordan UMR5208, B\^atiment Braconnier, 43 boulevard du 11 novembre 1918, 69622 Villeurbanne-cedex, France}
\email{wagner@math.univ-lyon1.fr}
\thanks{Recherche partiellement soutenu par le projet ANR-09-BLAN-0047 Modig.}
\date{\today}
\begin{abstract} Nous généralisation la construction de Frécon du {\em radical inévitable} aux groupes dans les théories stables et même simples.\par
{\noindent\sc Abstract.} We generalize Frécon's construction of the {\em inevitable radical} to groups in stable and even simple theories.\end{abstract} 
\maketitle

\section{Introduction}
Dans \cite{frecon08} Olivier Frécon a défini un sous-groupe définissablement caractéristique d'un groupe de rang de Morley fini, le sous-groupe $\In(G)$ des éléments {\em inévitables}, sous-groupe minimal tel qu'on puisse espérer construire une géométrie sur $G/\In(G)$. Nous allons généraliser ses définitions au cas d'un groupe stable, ou encore hyperdéfinissable dans une théorie simple, et étudier les propriétés des sous-groupes obtenus.

\section{Les définitions inévitables}
Rappelons d'abord les définitions (ou plutôt les caractérisations équi\-valentes) de \cite{frecon08}:
\defn Soit $G$ un groupe de rang de Morley fini. Une famille définissable $\F$ de sous-groupes connexes de $G$ est {\em géo\-métrique} si l'ensemble $\{g\in G:\exists!\ F\in\F\ g\in F\}$ est générique dans $G$.\par
Un élément $g\in G$ est {\em géométrique} s'il y a une famille géométrique $\F$ avec $g\notin\bigcup\F$.\par
$G$ est {\em géométrique} si tout élément non-trivial de $G$ est géométrique.\par
Un élément non-géométrique est {\em inévitable}. L'ensemble des éléments inévitables est noté $\In(G)$.\edefn
Frécon montre que $\In(G)$ est un sous-groupe définissable et définis\-sablement caractéristique de $G$. En plus, si $\F$ est géométrique pour $G$, alors $\{F\in\F:\In(G)\le F\}$ est géométrique, et donc $G/\In(G)$ est un groupe géométrique.
\bem Puisque la famille $\F$ est définissable, les sous-groupes dans $\F$ sont uniformément définissables. Par contre, comme on ne sait pas si la connexité est une propriété définissable, on ne sait pas non plus si toute famille de sous-groupes connexes uniformément définissables est contenue dans une famille définissable de sous-groupes connexes. Plus généralement on ne sait pas si la famille des composantes connexes d'une famille de sous-groupes uniformément définissables est encore uniformément définissable. Il n'est donc pas clair si $\In(G)$ est le plus petit sous-groupe normal de $G$ tel que le quotient soit géométrique.\ebem
Si $\F$ est une famille géométrique d'un groupe de rang de Morley $G$, alors on peut définir une géométrie naturelle dont les points sont les éléments de $G$, et les droites les translates (à gauche) des sous-groupes dans $\F$. Alors génériquement deux points de $G$ sont sur une droite commune unique.

Frécon note qu'on doit aussi considérer les familles géométriques des puissances cartésiennes de $G$~: Pour le groupe additif d'un pur corps $K$ algébriquement clos, $\In(K)=K$ comme pour tout groupe minimal, mais $\In(K\times K)$ est trivial.
\defn\label{defnInP} Soit $G$ un groupe de rang de Morley fini et $\rho_n:G\to\prod_{i<n} G$ le plongement $g\mapsto(g,1,\ldots,1)$.
On pose
$$\In_P(G)=\bigcap_{n>0}\rho_n^{-1}(\In(\prod_{i<n}G)).$$\edefn
Alors $\In_P(G)$ est définissable, définissablement caractéristique, con\-tenu dans $\In(G)$, et vérifie $\In_P(G/\In_P(G))=1$. De plus, $$\In_P(\prod_{i<n}G)=\prod_{i<n}\In_P(G).$$
\bsp\cite{jaligot06} Soit $G$ un groupe de rang de Morley fini et $C$ un sous-groupe de {\em Carter} (sous-groupe définissable connexe nilpotent presque auto-normalisant) {\em généreux} (dont la réunion des conjugués recouvre $G$ génériquement). Alors la famille des conjugués de $C$ est géométrique.\ebsp
\bsp\cite{frecon08} Soit $G$ un groupe de rang de Morley fini et $H$ un sous-groupe définissable sans torsion, d'indice fini dans son normalisateur. Alors soit $H$ admet une famille géométrique de sous-groupes propres, soit la famille de conjugués de $H$ est géométrique dans $G$.\par
En particulier, si $H=C$ est un sous-groupe de Carter et $G$ est un groupe simple connexe minimal, alors soit la famille de conjugués de $C$ est géométrique, soit $C$ a une famille géométrique de sous-groupes propres et interprète un corps algébriquement clos.\ebsp
Nous conseillons l'article de Frécon \cite{frecon08} pour une analyse des groupes géométriques algébriques et une discussion des groupes géométriques de rang de Morley fini. Pour plus de renseignements sur les groupes stables, auxquels nous cherchons à généraliser les résultats de Frécon, le lecteur pourra consulter \cite{poizat87}, \cite{wagner97} ou \cite{wagner00a}; quant aux théories simples et leurs groupes hyperdéfinissables il y a \cite{wagner00}, \cite{wagner01a} et \cite{wagner05}.

\section{Rajoutons un peu de simplicite !}
Les définitions ci-dessus restent raisonnables dans le contexte plus général des groupes $\omega$-stables, où les composantes connexes existent et sont définissables. Dans un groupe stable, les composantes connexes ne sont a priori que type-définissables, données par une intersection infinie de sous-groupes définissables d'indices finis~; de plus, le contexte naturel dans une théorie stable sont les groupes type-définisables. Pire encore, pour un groupe hyperdéfinissable dans une théorie simple, les composantes connexes dépendent des paramètres, et même dans le cas type-définis\-sable la question si un tel groupe est intersection de groupes définis\-sables reste ouverte.
\defn Soit $G$ un groupe hyperdéfinissable sur $A$, et $H$ un sous-groupe de $G$ hyperdéfinissable sur $B$.\begin{itemize}
\item Un sous-groupe hyperdéfinissable $K$ de $G$ est {\em commensurable} avec $H$ si $H\cap K$ est d'indice borné dans $H$ et dans $K$.
\item $H$ est {\em localement connexe} si $H$ est égal à tout $H^\gamma$, conjugué de $H$ par un élément de $G$ ou par un automorphisme modèle-théorique fixant $A$, dès que $H$ et $H^\gamma$ sont commensurables.
\item $B$ est le {\em paramètre canonique} de $H$ (sur $A$) si tout automor\-phis\-me fixant $A$ stabilise $H$ si et seulement s'il fixe $B$.
\end{itemize}\edefn
La notion de connexité locale dépend du groupe ambient $G$ et de ses paramètres $A$.
\tats\begin{itemize}
\item Soit $G$ un groupe type-définissable sur $A$ dans une théorie stable, et $H\le G$ un sous-groupe type-définissable. Alors $H$ possède un paramètre canonique, et sa composante connexe $H^0$ est localement connexe. Si $H$ est relativement définissable, alors $H$ possède un paramètre canonique imaginaire fini, et l'intersection de tous les $A$-conjugués et de tous les $G$-con\-ju\-gues de $H$ commensurables avec $H$ est un sous-groupe $H^{lc}$ relativement défi\-nis\-sable localement connexe d'indice fini dans $H$.
\item\cite[Corollary 4.2.10, Corollary 4.5.16 et Lemma 4.5.19]{wagner00} Soit $G$ un groupe hyperdéfinissable sur $A$ dans une théorie simple, et $H\le G$ un sous-groupe hyperdéfinissable. Alors il existe un sous-groupe $H^{lc}$ hyperdéfinissable localement connexe commensurable avec $H$~; si $G$ est type-définissable et $H$ relativement définissable, alors $H^{lc}$ est relativement définissable. Tout groupe localement connexe à un paramètre canonique.
\end{itemize}\etats
\bem Soit $F=\{G\cap F_c:c\models p\}$ une famille de sous-groupes relativement définissables localement connexes dans une théorie simple telle que $F_c$ définisse un groupe pour tout $c$ (par exemple si le groupe ambient $G$ est définissable). Pour $c\models p$ soit $n(\p,k)=D^*(G\cap F_c,\p,k)$ le $(\p,k)$-rang local de $G\cap F_c$. Alors si $c', c''\models p$ et 
$$D^*(G\cap F_{c'}\cap F_{c''},\p,k)\ge n(\p,k)$$
pour tout $(\p,k)$, alors $G\cap F_{c'}$ et $G\cap F_{c''}$ sont commensurables, et donc $c'=c''$. Par compacité il y a un ensemble fini $\{(\p_i,k_i):i<n\}$, une formule $\psi_0\in p$ et un ensemble définissable $X\supseteq G$ tels que
$$\psi_0(c')\land\psi_0(c'')\land\bigwedge_{i<n}D^*(X\cap F_{c'}\cap F_{c''},\p_i,k_i)\ge n(\p_i,k_i)$$
implique $c'=c''$. Soit $\psi(x)$ la formule
$$\psi_0(x)\land\bigwedge_{i<n}D^*(X\cap F_x,\p_i,k_i)\ge n(\p_i,k_i).$$
Alors la  famille $\bar\F=\{G\cap F_c:c\models\psi\}$ est une famille définissable de sous-groupes localement connexes de $G$.\ebem
A partir de maintenant soit $G$ un groupe hyperdéfinissable sur $\emptyset$ dans une théorie simple.
Comme nos familles de sous-groupes ne consistent plus de groups connexes, mais localement connexes, il convient de relacher un peu la condition d'unicité dans la définition d'une famille géométrique.
\defn\label{defgeom} Soit $A$ un ensemble de paramètres. La famille $\F$ des $A$-conjugués d'un sous-groupe localement connexe hyperdéfinissable de $G$ est {\em (pseudo-)géo\-mé\-trique} si $\bigcup\F$ est générique, et si pour tout $g\in G$ générique l'ensemble $\{F\in\F:g\in F\}$ est de cardinal au plus $1$ (resp.\ de cardinal borné).\par
Un élément $g\in G$ est {\em (pseudo-)géométrique} s'il y a une famille (pseudo-)géométrique $\F_g$ avec $g\notin\bigcup\F_g$. Le groupe $G$ est {\em (pseudo-) géométrique} si tout élément non-trivial de $G$ est (pseudo-)géométrique.\par
Un élément non-(pseudo-)géométrique est {\em (pseudo-)inévitable}.\par
L'ensemble des éléments inévitables est noté $\In(G)$, l'ensemble des éléments pseudo-inévitables $\Pin(G)$.\edefn
Si $\F$ est une famille (pseudo-)géométrique, puisque tous ses éléments sont conjugués, la connexité locale implique que $F_1$ et $F_2$ dans $\F$ sont égaux dès qu'ils sont commensurables. On notera aussi que tout membre de $\F$ contient un élément générique de $G$ (qui algébraise son paramètre canonique, bien sur), comme tous les groupes dans $\F$ sont conjugués par $A$-automorphisme.
\bem Soit $\F$ une famille hyperdéfinissable sur $A$ avec $\bigcup\F$ générique et $\{F\in\F:g\in F\}$ de cardinal au plus un (resp.\ borné) pour tout $g\in G$ générique sur $A$. Soit $F_0\in\F$ et $g\in F_0$ générique de $G$ sur $A$. Soit $\F_0\subseteq\F$ la sous-famille des $A$-conjugués de $F_0$. Alors $\F_0$ est une famille (pseudo-)géométrique pour $G$. La condition que les éléments de $\F$ soient conjugués sur $A$ n'est donc pas une vraie restriction.\ebem
\bem Soit $\F=\{G\cap F_c:c\models p\}$ une famille géométrique de sous-groupes relativement définissables dans une théorie simple telle que $F_c$ définisse un groupe pour tout $c$. Alors par compacité il existe un ensemble générique $X$ de $G$ et une formule $\p\in p$ tel que pour tout $c'\models\p$ l'ensemble $G\cap F_{c'}$ est un sous-groupe localement connexe de $G$, et pour tout $g\in X$ il existe un unique $c_g\models\p$ avec $g\in G\cap F_{c_g}$. Si donc $G$ est définissable, alors $\{G\cap F_{c'}:c'\models\p\}$ est une famille géométrique au sens de Frécon (sauf qu'on a remplacé la connexité par la connexité locale, faute de pouvoir obtenir cette première définissablement). Réciproquement, si $\F$ est une famille géo\-métrique sur $A$ au sens de Frécon, $g\in G$ est générique sur $A$ et $F_c\in\F$ contient $g$, alors la famille des $A$-conjugués de $F_c$ est une famille géomérique au sens de la définition \ref{defgeom}.\ebem
\bem Si $\F$ consiste de sous-groupes connexes, alors l'élément $g\in G$ générique qui est contenu dans un membre de $\F$ est forcément dans la composante connexe de $G$~; si la théorie est stable, il n'y a qu'un seul type possible, et on peut supposer dès le départ que $G$ soit connexe. En particulier notre définition \ref{defgeom} généralise bien celle de Frécon.

Par contre, même dans le cas stable, si $\F$ ne consiste pas de sous-groupes connexes, il n'y a pas raison que $g$ soit générique principal, ni même qu'un générique principal soit contenu dans un membre de $\F$~: Soit $G$ le produit semidirect de $\QQ$ avec une involution $i$ qui agit par inversion, et $\F$ la famille de conjugués de $\langle i\rangle$. Alors $G$ est $\omega$-stable de rang $1$ et degré $2$~; un générique non-principal est dans un seul membre de $\F$, et $\bigcup\F\cap G^0=\{0\}$.\ebem
\bem Afin d'éviter des cas dégénérés on peut demander dans la définition d'une famille géométrique $\F$ que les groupes $F\in\F$ soient infinis, ou qu'un générique de $G$ ne soit pas étranger à $F$, ou encore que $G$ soit $F$-analysable (conditions de plus en plus restrictives). Inversement, nos preuves ne se servent de la connexité locale que pour s'assurer qu'un paramètre canonique existe~; en utilisant des paramètres de définition quelconques raisonnablement indépendants on peut se passer de cette hypothèse.\ebem
\bem Si $\F$ est une famille géométrique pour $G$, alors la famille des translatés à gauche des sous-groupes dans $\F$ nous donne une géométrie sur $G$ tel que génériquement deux points sont sur au plus une droite commune. Si $\F$ n'est que pseudo-géométrique, cela nous donne une pseudo-géométrie~: Génériquement deux points ne sont que sur un nombre borné de droites communes. Par compacité, si le paramètre canonique d'un groupe dans $\F$ est un élément imaginaire (par exemple si les groupes sont relativement définissables), ce nombre borné est fini.\ebem
\bem Si $\F_i$ est une famille (pseudo-)géométrique de $G_i$ pour $i=1,2$ et $F_1\in\F_1$ et $F_2\in\F_2$ ont des paramètres indépendants, alors la famille des conjugués de $F_1\times F_2$ est (pseudo-)géométrique dans $G_1\times G_2$.\ebem
\satzli\label{extfamille} Soit $A\subseteq B$, et $F\le G$ un sous-groupe localement connexe hyperdéfinissable de paramètre canonique $c$ avec $c\ind_AB$. Si la famille $\F_A$ des $A$-conjugués de $F$ est (pseudo-)géométrique, alors la famille $\F_B$ des $B$-conjugués l'est. La réciproque est vrai pour les familles pseudo-géométriques~; dans le cas d'une famille géométrique il faut supposer en plus que $\tp(c/A)$ soit stationnaire (par exemple si $A=\acl^{eq}(A)$ dans une théorie stable).\esatzli
\bew Soit $\F_A$ (pseudo-)géométrique, et $g\in F$ générique de $G$ sur $A$. On peut supposer $g\ind_{Ac}B$. Par transitivité $B\ind_A gc$, et $g$ est générique sur $B$. Donc $\bigcup\F_B$ est générique. Mais tout générique de $G$ sur $B$ est générique sur $A$, et $\F_B\subseteq\F_A$. Donc $\F_B$ est \hbox{(pseudo-)} géométrique.

Pour la réciproque, supposons $\F_B$ pseudo-géométrique. Evidemment $\bigcup\F_A\supseteq\bigcup\F_B$ est générique dans $G$. Soit $g\in F$ générique de $G$ sur $A$~; on peut supposer $g\ind_A B$. Soient $c_i\equiv_Ac$ les paramètres canoniques des $F_i\in\F_A$ avec $g\in F_i$ pour $i\in I$. Si $I$ n'est pas borné, on peut supposer $(c_i:i\in I)$ indiscernable sur $A$ et $(c_i:i\in I)\ind_{Ag}B$, d'où $B\ind_A (c_i:i\in I)$. Comme $p(X,c)=\tp(B/Ac)$ ne devie pas sur $A$, on trouve une réalisation 
$$B'\models\bigcup_{i\in I}p(X,c_i)\quad\mbox{avec}\quad B'\ind_A (c_i:i\in I)\ ;$$
on peut supposer en plus que
$$B'\ind_{(A,c_i:i\in I)}g,$$
d'où $B'\ind_A g$. Donc $g$ est générique sur $B'$ et $F_i\in\F_{B'}$, l'image de $\F_B$ sous un $A$-automorphime envoyant $B$ sur $B'$. Comme $\F_B$ est pseudo-géométrique, $\F_{B'}$ l'est aussi, est $I$ doit être bornée après tout.

Soit enfin $\F_B$ géométrique et $\tp(c/A)$ stationnaire. Alors $\bigcup\F_A$ est générique~; soient donc $g\in F$ et $(F_i,c_i:i\in I)$ comme dans le paragraphe précédent. On suppose encore
$$g\ind_A B\quad\mbox{et}\quad(c_i:i\in I)\ind_{Ag}B.$$
Donc $c_i\ind_AB$ pour tout $i\in I$~; par stationarité $c_i\equiv_B c$ et $F_i\in\F_B$. Puisque $g$ est générique sur $B$ et $\F_B$ est géométrique, $|I|=1$. Donc $\F_A$ est géométrique.\qed

\bem En particulier, dans une théorie simple $\In(G)$ et $\Pin(G)$ sont invariants sous l'adjonction de paramètres.\ebem

\kor\label{korindep} Soit $\F$ une famille (pseudo-)géométrique pour $G$. Si $g\in G$ est (pseudo-)inévitable, alors $g\in F$ pour tout $F\in\F$ dont le paramètre canonique $c$ est indépendant de $g$ sur $A$.\ekor
\bew La famille $\F_g$ des $Ag$-conjugués de $F$ est toujours (pseudo-)géométrique, donc $g\in\bigcup\F_g$. Mais alors $g\in F$.\qed

\kor $\In(G)$ et $\Pin(G)$ sont des sous-groupes de $G$ définis\-sable\-ment et modèle-théoriquement caractéristiques (c'est-à-dire invariant par automorphisme définissable ou modèle-théorique).\ekor
\bew Si $\F$ est une famille géométrique pour $G$ et $g,g'\in\In(G)$, alors pour $F\in\F$ de paramètre canonique indépendant de $g,g'$ on a $g,g'\in F$, et donc $g'g^{-1}\in F$, d'ou $g'g^{-1}\in\In(G)$~: Il s'agit bien d'un sous-groupe. Puisque $\In(G)$ est invariant sous l'adjonction de paramètres né\-ces\-saires pour définir un automorphisme $\sigma$, et comme $\F$ est géométrique si et seulement si $\sigma(\F)$ est géométrique pour tout automorphisme définissable ou modèle-théorique, $\In(G)$ est définissablement et modèle-théoriquement caractéristique.

La preuve dans le cas pseudo-géométrique est analogue.\qed

\kor S'il y a une famille géométrique non-triviale, $\In(G)$ ne contient aucun élément générique. S'il y a une famille pseudo-géomé\-trique non-bornée, $\Pin(G)$ ne contient aucun élément générique.\ekor
\bew Supposons que $g\in\In(G)$ soit générique pour $G$, et soit $\F$ une famille géométrique sur $A$. On peut supposer que $A\ind g$. Soient $F, F'\in\F$ de paramètres canoniques $c,c'$. On peut les choisir tels que $c,c'\ind_A g$. Alors $g\in F$ et $g\in F'$ par le corollaire \ref{korindep}, d'où $F=F'$ et $\F$ est triviale.

Si $g\in\Pin(G)$ est générique et $\F$ est une famille pseudo-géométrique sur $A$, avec $A\ind g$, on considère des groupes distincts $F_i\in\F$ de paramètres canoniques $c_i$, pour $i\in I$. On peut supposer $(c_i)_{i\in I}\ind_A g$. Alors $g\in F_i$ pour tout $i\in I$ par le corollaire \ref{korindep}. Ainsi $I$ et donc $\F$ sont bornés.\qed

\satz Si $T$ est stable, $\In(G)$ et $\Pin(G)$ sont hyperdéfinissables sur $\emptyset$.\esatz
\bew  Soit $\F$ une famille géométrique, disons la famille des $A$-conjugués d'un groupe $F_c$ localement connexe de paramètre canonique $c$. On considère l'ensemble hyperdéfinissable
$$G_{\F}=\{g\in G:\exists\,y\models\stp(c/A)\ [\,y\ind_Ag\land g\in F_y\,].$$
Alors $G_\F$ est un sous-groupe hyperdéfinissable de $G$ contenu dans $\bigcup\F$.
Si $g\in G$ est inévitable, alors pour $c'\models\stp(c/A)$ avec $c'\ind_Ag$ le corollaire \ref{korindep} montre que $g\in F_{c'}$, d'où $g\in G_\F$ et $In(G)\subseteq G_\F$. Si $g\in G$ est géométrique, il existe une famille géométrique avec $g\notin\bigcup\F$, d'où $g\notin G_\F$. Ainsi $\In(G)=\bigcap_\F G_\F$, ce qui est hyperdéfinissable par la condition de chaîne dans les théories stables. Comme $\In(G)$ est $\emptyset$-invariant, il est hyperdéfinissable sur $\emptyset$.

La preuve pour $\Pin(G)$ est analogue.\qed

\kor Soit $T$ stable. Si $\F$ est une famille géométrique, alors $\In(G)\le\bigcap\F$~; si $\F$ est pseudo-géométrique, alors $\Pin(G)\le\bigcap\F$.\ekor
\bew Si $\F$ est $B$-invariant, $F\in\F$ avec paramètre canonique $c$, et $g\in In(G)$ est générique sur $Bc$, alors $g\ind_B c$, d'où $g\in F$. Donc $\In(G)\le F$, et $\In(G)\le\bigcap\F$.

La preuve dans le cas pseudo-géométrique est analogue.\qed

\kor Soit $T$ stable.\begin{itemize}
\item $\In(G)=\bigcap\,\{\bigcap\F:\F\mbox{ géométrique}\}$.
\item $\Pin(G)=\bigcap\,\{\bigcap\F:\F\mbox{ pseudo-géométrique}\}$.\qed\end{itemize}\ekor

\kor $\In$ et $\Pin$ sont des radicaux, c'est-à-dire $\In(G/\In(G))$ et $\Pin(G/\Pin(G))$ sont triviaux.\ekor
\bew Soit $\F$ géométrique pour $G$. Comme $\In(G)\le F$ pour tout $F\in\F$, l'image $\bar\F$ de $\F$ dans $G/\In(G)$ est géométrique. Comme $\bigcap_{\F}\bigcap\F=\In(G)$, on a
$$\In(G/\In(G))\le\bigcap_{\F}\bigcap\bar\F=\big(\bigcap_{\F}\bigcap\F\big)/\In(G)=\In(G)/\In(G)=\bar1.$$
La preuve pour $\Pin$ est analogue.\qed

\bem On peut alors définir $\In_P$ et $\Pin_P$ de la mème manière que Frécon (définition \ref{defnInP})
et obtenir que $\In_P(G)$ et $\Pin_P(G)$ sont des radicaux hyperdéfinissables, définissablement et modèle-théoriquement caractéristiques, et contenu dans $\In(G)$ et $\Pin(G)$, respectivement. De plus, $$\In_P(\prod_{i<n}G)=\prod_{i<n}\In_P(G)\qquad\mbox{et}\qquad\Pin_P(\prod_{i<n}G)=\prod_{i<n}\Pin_P(G).$$
\ebem

\section{Dans le cas simple, les choses se compliquent...}
Si $G$ est un groupe dans une théorie simple où on n'a la condition de chaîne qu'à indice borné près, il convient de remplacer les groupes $\In(G)$ et $\Pin(G)$ par des approximations $\aIn(G)$ et $\aPin(G)$, comme c'est déjà le cas pour le centralisateur, le normalisateur ou le centre approximatif \cite[Definition 4.4.9]{wagner00}.

\satz Il existe un sous-groupe normal $\aIn(G)$ hyper\-définis\-sable sur $\emptyset$ tel que $\In(G)\le\aIn(G)$, et si $\F$ est une famille géométrique, alors $F$ intersecte $\aIn(G)$ dans un sous-groupe d'indice borné pour tout $F\in\F$. De même, il existe un sous-groupe normal $\aPin(G)$ hyperdéfinis\-sable sur $\emptyset$ tel que $\Pin(G)\le\aPin(G)$, et si $\F$ est une famille pseudo-géométrique, alors $F$ intersecte $\aPin(G)$ dans un sous-groupe d'indice borné pour tout $F\in\F$.
Les groupes $\aIn(G)$ et $\aPin(G)$ sont invariants sous automorphisme $\emptyset$-définissable, et commensurables à leurs conjugués par un automorphisme définissable.\esatz
\bew Soit $F=F_c$ un groupe dans une famille géométrique, et $p_F$ le type Lascar-fort de son paramètre canonique $c$ sur $\emptyset$.
On considère l'ensemble hyperdéfinissable
$$X_F=\{g\in G:\exists\,y\models p_F\ [\,y\ind g\land g\in F_y\,]\}.$$
Soit $X=\bigcap_F X_F$. Puisque l'intersection porte sur les differentes possibilités pour $p_F\in S(\bdd(\emptyset))$, l'intersection est bornée et $X$ est hyperdéfinissable.

Si $g',g''\in X$ sont indépendants, alors pour tout $p_F$ il y a $c',c''\models p_F$ avec $c'\ind g'$ et $g'\in F_{c'}$, ainsi que $c''\ind g''$ et $g''\in F_{c''}$. Grâce au théorème d'indépendance on peut supposer $c'=c''$ et $c'\ind g,g'$. Alors $g,g'\in F_{c'}$ et $g'g^{-1}\in F_{c'}$~; comme $g'g^{-1}\ind c'$ on a $g'g^{-1}\in X_F$, d'où $g'g^{-1}\in X$. On pose $H=X^2$. Alors $H$ est un sous-groupe de $G$ hyperdéfinissable sur $\bdd(\emptyset)$ par \cite[Lemma 4.4.8]{wagner00}, et $X$ contient tous les types génériques de $H$.

Si $g\in G$ est inévitable et $\F_A$ est une famille géométrique sur $A=\bdd(A)$, on choisit $A'\equiv_{\bdd(\emptyset)}A$ avec $g\ind A'$. Alors $\F_{A'}$ est toujours géométrique~; si $F'\in\F_{A'}$ est de paramètre $c'$ canonique indépendant de $g$ sur $A'$, alors $g\in F'$ par le lemme \ref{korindep}. Comme $c'\ind g$ par transitivité, on a $g\in X_{F'}=X_F$, et donc $g\in X$, d'où $\In(G)\subseteq X\subseteq H$.

Soit $F\in\F_A$ de paramètre canonique $c$, et $h'\in H$ générique sur $Ac$. Alors il y a $c'\models p_F$ avec $h'\in F_{c'}$ et $c'\ind h'$. Soit $A'\ind_{c'}h'$ tel que $A'c'\equiv_{\bdd(\emptyset)}Ac$, et $\sigma$ un $\bdd(\emptyset)$-automorphisme qui envoit $A'c'$ sur $Ac$. Soit $h=\sigma(h')$. Comme $h'\ind A'c'$ on a $h\ind Ac$, et $h$ est générique dans $H$ sur $Ac$. Comme $h\in F$, ceci signifie que $F$ intersecte $H$ dans un sous-groupe d'indice borné.

A priori la définition de $H$ dépend des paramètres~; on dénotera par $H_B$ le groupe obtenu par la même construction mais avec des paramètres $B$ nommés. On montrera que $H_B$ est un sous-groupe d'indice borné dans $H$. Soit donc $p_F$ le type Lascar-fort sur $\emptyset$ du paramètre canonique $c$ d'un groupe $F_c$ dans une famille géométrique. Soit $p'_F$ une extension non-deviante de $p_F$ sur $\bdd(B)$~; par la proposition \ref{extfamille} c'est le type Lascar-fort du paramètre canonique d'une famille géométrique dont les paramètre incluent $B$. Alors si $h\in H_B$ est générique sur $B$, il y a $c'\models p'_F$ avec $h\in F_{c'}$ et $c'\ind_B h$. Mais $c'\ind B$, d'où $c'\ind h$ et $h\in X\subseteq H$. Donc $H_B\le H$.

Réciproquement, pour tout $p'_F$ type Lascar-fort sur $B$ du paramètre canonique d'un groupe dans une famille géométrique (avec $B$ nommé) soit $c_F\models p'_F$ une réalisation. Alors $F_{c_F}$ intersecte $H$ dans un sous-groupe d'indice borné, et il y a $h\in H$ générique sur $B\cup(c_F)_F$ avec $h\in F_{c_F}$ pour tout $F$. Donc $h\ind B\cup(c_F)_F$, d'où $h\ind_B c_F$ pour tout $F$, et $h\in X_B\subseteq H_B$. Ceci signifie que $H_B$ est générique, et donc d'indice borné, dans $H$.

Puisque l'image d'une famille géométrique sous un automorphisme $\bdd(\emptyset)$-définissable est encore géométrique avec les mêmes paramètres canoniques, $H$ est invariant sous automorphisme $\bdd(\emptyset)$-définissable. Si donc un automorphisme $\sigma$ est définissable à l'aide de paramètres $B$, on a que $H_B=\sigma(H_B)$ est d'indice borné dans $H$ et dans $\sigma(H)$. Ainsi $H$ et $\sigma(H)$ sont commensurables.

Soit $p$ un type générique principal Lascar-fort de $G$. On pose
$$Y=\{g\in G:\exists\,y\models p\ [\,y\ind g\land g\in H^y\,]\}.$$
Comme avant, $Y^2=K$ est un sous-groupe de $G$ et $Y$ contient tous les génériques de $K$. Comme $H$ contient $\In(G)$ et ce dernier est normal, $\In(G)\le H^y$ pour tout $y\models p$, et $\In(G)\le K$. Si $y\models p$, alors $H\cap H^y$ est d'indice borné dans $H$ et contient un générique $h$ de $H$. Alors $h\ind y$ et $h\in K$. Ainsi $K$ intersecte $H$ dans un sous-groupe d'indice borné. Réciproquement, si $h\in Y$ est générique dans $K$, il y a $y\models p$ avec $y\ind h$ et $h\in H^y$. Mais alors $\tp(h)$, $\tp(h/y)$ et $\tp(h^{y^{-1}}/y)$ ont les mêmes rangs locaux stratifiés~; comme $h^{y^{-1}}\in H$ on conclut que $H\cap K$ est d'indice borné dans $K$~: On a bien que $H$ et $K$ sont commensurables.

Soit $g$ générique principal de $G$, et $h\in K$ générique sur $g$. Alors il y a $y\models p$ avec $h\in H^y$ et $y\ind h$. Puisque $p$ est générique et $g$ générique principal, il y a $y'\models p$ avec $y'g\models p$ et $y'\ind g$. Comme $h\ind g$ on peut prendre $y=y'\ind h,g$ par le théorème d'indépendance. Alors $h^g\in H^{yg}$ avec $yg\models p$ et $yg\ind_g h^g$~; puisque $yg\ind g$ par généricité on a $yg\ind h^g$, d'où $h^g\in Y\subseteq K$. Ainsi $K$ est normalise par tous les génériques principaux de $G$, et donc par $G^0$.

Soit $\aIn(G)$ l'intersection de $K$ avec tous ses $G$-conjugués, et tous ses conjugués par $\emptyset$-automorphisme. C'est une intersection bornée~; $\aIn(G)$ est donc hyperdéfinissable sur $\emptyset$, et clairement normal dans $G$. Comme tous ces conjugués sont commensurables, $\aIn(G)$ est commensurable avec $K$, et donc avec $H$. Enfin, $\aIn(G)$ est commensurable avec tout conjugué par automorphisme définissable, et tout groupe $F$ dans une famille géométrique l'intersecte dans un sous-groupe d'indice borné, puisque c'est vrai de $H$.

La définition de $\aPin(G)$ et la preuve de ses propriétés sont analogues, en utilisant les familles pseudo-géométriques.\qed

\bem Il se peut que $\In(G)$ soit trivial, bien que $\aIn(G)$ est non-borné~; de plus, l'effet de la saturation n'est pas évident, puisqu'un groupe non-saturé a certes moins d'éléments dans ses sous-groupes hyperdéfinissables, mais aussi moins de famille géométriques qui interdiraient aux éléments d'appartenir à $\In(G)$. En particulier, si $G\preceq G^*$, il n'est pas clair si $\In(G)=\In(G^*)\cap G$, ni même si $\In(G)\le\In(G^*)$.\ebem

\satzli S'il y a une famille géométrique non-triviale, $\aIn(G)$ est d'indice non-borné dans $G$. S'il y a une famille pseudo-géométrique non-bornée, $\aPin(G)$ est d'indice non-borné dans $G$.\esatzli
\bew Supposons $\aIn(G)$ d'indice borné dans $G$. Soit $\F$ une famille géométrique et $F, F'\in\F$. Alors $\aIn(G)\cap F\cap F'$ est d'indice borné dans $G$, et contient un générique $g$. Mais $g\in F$ et $g\in F'$, d'où $F=F'$.

Si $\aPin(G)$ est d'indice borné dans $G$ et $\F$ est une famille pseudo-géométrique, soient $(F_i:i\in I)$ des groupes distincts dans $F$. Alors $\aPin(G)\cap\bigcap_{i\in I}F_i$ est générique dans $G$ et contient un générique $g$. Donc $I$ et $\F$ sont bornés, puisque $g\in F_i$ pour tout $i\in I$.\qed

\bem Si $\F$ est une famille géométrique et $F\in\F$ coupe $\aIn(G)$ dans un sous-groupe propre, soit $g\in F$ générique pour $G$ sur $\emptyset$ et $h\in\aIn(G)\setminus F$ avec $g\ind h$. Alors $gh\notin F$; mais $gh$ est toujours générique, et peut être contenu dans un autre groupe $F'\in\F$. Notons que si $F\cdot\aIn(G)=F'\cdot\aIn(G)$, alors $F$ et $F'$ sont commensurables et donc égaux. En particulier la famille $\bar\F$ des images de $\F$ dans $G/\aIn(G)$ n'est plus géométrique, puisque le générique $g\cdot\aIn(G)$ est contenu dans plusieurs groupes dans $\bar\F$.\ebem

\satzli Si $\F$ est une famille pseudo-géométrique dans $G$, alors l'image $\bar F$ de $\F$ dans $G/\aPin(G)$ est une famille pseudo-géométrique. En particulier $\aPin$ est un radical~: $\aPin(G/\aPin(G))$ est trivial.\esatzli
\bew Comme $F$ et $F\cdot\aPin(G)$ sont commensurables pour $F\in\F$, si $F\cdot\aIn(G)$ et un conjugué groupe- ou modèle-théorique $(F\cdot\aIn(G))^\sigma$ sont commensurables, alors $F$ et $F^\sigma$ sont commensurables et donc égaux par connexité locale. Puisque $\aIn(G)$ est $\sigma$-invariant, $\bar\F$ consiste de groupes localement connexes. Pour montrer que $\bar\F$ est pseudo-géométrique, il suffit alors à montrer que tout générique de $G/\aPin(G)$ n'est contenu que dans un nombre borné de groupes dans $\bar\F$.

Soit $g\in G$ générique sur les paramètres $A$ qui servent à définir $\F$, et $H_g$ l'intersection avec $\aPin(G)$ de tous les $F\in\F$ qui contiennent $g$. C'est une intersection bornée~; $H_g$ est donc d'indice borné dans $\aPin(G)$. Soit $p$ un type complèt sur $A,g$ tel que $p(x)$ implique que $x\in\aPin(G)$ et $gx$ est générique dans $G$ sur $A$. Si $a,b\models p$ et $a\in bH_{gb}$, alors $ga$ est dans tous les $F\in\F$ qui contiennent $gb$~; par conséquent $H_{ga}\le H_{gb}$ et
$$aH_{ga}\subseteq aH_{gb}=bH_{gb}.$$
Montrons l'égalité. Sinon $x\in yH_{gy}$ définit un ordre partiel $\le$ sur $p$ avec $a<b$. Puisque $p$ est complet, un trouve une chaîne $a_0<a_1<a_2<\cdots$ de réalisations de $p$~; par compacité il y a $a_\omega\models p$ avec $a_\omega\ge a_i$ pour tout $i<\omega$, et puisque la chaîne $(a_i:i<\omega)$ est strictement croissante, $a_\omega>a_i$ pour tout $i<\omega$. Récursivement on trouve ainsi une chaine stricte aussi longue qu'on veut, et par le théorème d'Erdös-Rado il y a une sous-chaîne $4$-indiscernable infinie, qu'on nomme encore $(a_i:i\le\omega)$. Soit $\p(x,yz)$ une formule dans $y\le x\le z$ tel que
$$\tp(a_0a_1a_2a_3)\land p(x)\land p(x')\land\p(x,a_0a_1)\land\p(x',a_2a_3)$$
implique $x\not=x'$. Alors
$$\begin{aligned}D(p(x)\land a_0\le x\le a_1,\p,2)&=D(p(x)\land a_0\le x\le a_\omega,\p,2)\\
&\ge D(p(x)\land a_0\le x\le a_1,\p,2)+1,\end{aligned}$$
une contradiction. Donc $H_{ga}=H_{gb}$ et $aH_{ga}=bH_{gb}$. Ainis $ga$ et $gb$ sont contenus dans les mêmes groupes $F\in\F$.

Autrement dit, pour $h\models p$ les ensembles $hH_{gh}$ définissent une partition de $p$~; soit $\sim_p$ la relation d'équivalence correspondante. Mais comme $H_{gh}$ est d'indice borné dans $\aPin(G)$, pour $h\models p$ la classe $hH_{gh}$ contient un élément générique $h'$ de $\aPin(G)$ sur $A,g,h$. Donc $h'\ind_{A,g}h$, et $h\ind_{A,g}h_{\sim_p}$ puisque $h_{\sim_p}=h'_{\sim_p}$. Ainsi $h_{\sim_p}\in\bdd(A,g)$ et il n'y a qu'un nombre borné de $\sim_p$-classes.

Or, pour tout $h\in\aPin(G)$ et $F\in\F$ tel que $gh\in F$, comme $F$ coupe $\aPin(G)$ dans un sous-groupe d'indice fini, il y a $h'\in\aPin(G)\cap F$ générique sur $A,g,h$. Alors $hh'$ est générique dans $\aPin(G)$ sur $A,g$, et $ghh'\in F$ est générique dans $G$ sur $A$. Comme le nombre de types $p$ sur $A,g$ est borné, et chaque type ne contient qu'un nombre borné de $\sim_p$-classes, et chaque $\sim_p$-classe consiste d'éléments $h\models p$ tel que $gh$ est dans les mêmes groupes $F\in\F$, le translaté $g\aPin(G)$ n'intersecte qu'un nombre borné de groupes $F\in\F$. Puisque tout générique de $G/\aPin(G)$ se relève en un générique de $G$, la famille $\bar\F$ est pseudo-géométrique.

Enfin, un élément $\bar g\in G/\aPin(G)$ est dans un groupe $\bar F_c$ d'une famille pseudo-géométrique $\bar\F$ pour un $c\ind\bar g$ si et seulement si un pré-image $g\in G$ est dans le pré-imagé $F_c\cdot\aPin(G)$ de la famille $\F$ avec $c\ind g$. Alors $F_cg\cap\aPin(G)$ est non-vide et ainsi d'indice borné dans $\aPin(G)$~; soit $g'\in F_c g\cap\aPin(G)$ générique dans $\aPin(G)$ sur $c,g$. Alors $g'\ind g,c$, d'où $g'\ind_g c$ et $c\ind g,g'$. Ainsi $gg^{\prime-1}\ind c$ et $gg^{\prime-1}\in F_c$. Donc $gg^{\prime-1}\in X_F$, l'ensemble utilise dans la construction de $\aPin(G)$. On peut faire la même chose simultanément pour tous les $p_F$ possibles~: Il suffit de prendre $g'$ générique dans
$$\aPin(G)\cap\bigcap_{F}F_{c_F}g$$
où l'intersection $\bigcap_{F}$ porte sur un ensemble de représentants des groupes $F$ pour les types Lascar-forts possibles de leur paramètre canonique $c_F$, et $c_F\ind g$ avec $g\in F_{c_F}\cdot\aPin(G)$. Ainsi, $gg^{\prime-1}\in X$ et
$$\aPin(G/\aPin(G))\le X^2/\aPin(G).$$
Or, $X^2/\aPin(G)$ est un groupe borné. On a donc $Y/\aPin(G)=X^2/\aPin(G)$. Comme $\aPin(G)$ est l'intersection de tous les conjugues de $Y^2$ groupe- ou modèle-théoriques, l'intersection de tous les conjugues groupe- ou modèle-théoriques de $Y^2/\aPin(G)$ est bien triviale.\qed

\bem On peut enfin définir $\aIn_P$ et $\aPin_P$ de la mème manière que Frécon (définition \ref{defnInP}). Ce seront des sous-groupes hyperdéfinis\-sables normaux, et contenus dans $\aIn(G)$ et $\aPin(G)$, respectivement. De plus, $$\aIn_P(\prod_{i<n}G)=\prod_{i<n}\aIn_P(G)\qquad\mbox{et}\qquad\aPin_P(\prod_{i<n}G)=\prod_{i<n}\aPin_P(G).$$
\ebem

\end{document}